\documentclass{amsart}
\usepackage{amssymb,amsmath,latexsym}
\newtheorem{theorem}{Theorem}

\title{Quantized open chaotic systems}

\author{S. Nonnenmacher}

\address{Institut de Physique Th\'eorique\\
CEA/DSM/PhT, Unit\'e de recherche associ\'ee au CNRS\\
CEA-Saclay\\
91191 Gif-sur-Yvette, France}
\email{snonnenmacher@cea.fr}

\begin{document}

\maketitle

\begin{abstract}
Two different ``wave chaotic'' systems, involving complex eigenvalues
or resonances, can be 
analyzed using common semiclassical methods. In particular, one obtains fractal Weyl
upper bounds for the density of resonances/eigenvalues near the real axis, and a
classical dynamical criterion for a spectral gap.
\end{abstract}


\section{Introduction}
These notes present a sketch of semiclassical methods which can be
used to describe the  spectral properties (and as a consequence, the
long time properties) of a certain class of 1-particle
quantum chaotic systems. Here are two examples:
\begin{itemize}
\item damped waves $\psi(x,t)$ on a compact riemannian manifold $X$ of negative
  sectional curvature. The dynamics is described by
  the damped wave equation
\begin{equation}\label{e:damped}
\big(\partial^2_t -\Delta_X + 2b(x)\partial_t\big) \psi(x,t)=0\,,
\end{equation}
and the {\it damping function} $b(x)\geq 0$ is assumed to be
smooth. 
\item quantum scattering on $\mathbb{R}^d$, described
  by the Schr\"odinger equation 
\begin{equation}\label{e:Schro-t}
i\hbar\partial_t \psi(x,t)= P(\hbar) \psi(x,t)\,,\quad
P(\hbar)=-\frac{\hbar^2\Delta}2+V(x)\,,
\end{equation}
where the potential $V(x)$ has
compact support and consists of $3$ peaks centered on an
equilateral triangle; $\hbar$ is Planck's constant.
\end{itemize}
These two systems seem very different. The configuration
spaces are respectively compact and of infinite volume, the wave
equation does not depend on $\hbar$; the damped wave equation
can be rewritten in terms of a contracting semigroup, while the
propagator $e^{-itP(\hbar)/\hbar}$ is unitary.

In the damped wave situation, Eq.~\eqref{e:damped} can be diagonalized by a discrete set of
metastable modes $e^{-ik_j t}\psi_j(x)$, where 
$\psi_j(x)\in L^2(X)$ satisfies the generalized
eigenvalue equation
\begin{equation}\label{e:k_j}
\big(\Delta + k_j^2 + 2 i\,b(x)\,k_j\big) \psi_j(x)=0\,.
\end{equation}
The eigenvalues $k_j$ are complex, and lie in the strip $-2\max(b)\leq
\Im k_j \leq 0$. We are interested in the distribution of
these eigenvalues in the {\it high-frequency limit} $\Re k_j\to\infty$. 

In the scattering situation \eqref{e:Schro-t}, the
Hamiltonian $P(\hbar)$ is selfadjoint on $L^2(\mathbb{R^d})$, with
absolutely continuous spectrum on $\mathbb{R}_+$. However, the
Green's function $(P(\hbar)-z)^{-1}(x,y)$, well-defined for $\Im z>0$,
admits a meromorphic continuation through
$\mathbb{R}_+$ to the lower half-plane, with discrete poles $\{z_j(\hbar)\}$ of finite multiplicities (the
{\it resonances} of $P(\hbar)$). We will investigate the distribution of
these resonances, in the vicinity of a fixed energy $E>0$, in
the {\it semiclassical limit} $\hbar\to 0$. To each resonance is
associated a metastable state $\psi_j(x)$, which is not in $L^2$
but formally decays with time as 
$e^{-itz_j(\hbar)/\hbar}\psi_j(x)$.

The {\it lifetime} of a
  metastable state is given by the {\it imaginary part} of the
  eigenvalue, $\tau_j = \frac{1}{2|\Im k_j|}$, resp. $\tau_j =
  \frac{\hbar}{2|\Im z_j(\hbar)|}$.

Here are some common features of the two systems.
In \S\ref{s:damped1} we show that the high-frequency limit $\Re k_j\to\infty$ is similar with
  the semiclassical limit $\hbar\to 0$. It is then relevant to study the corresponding {\it classical dynamics}:
  in the damped wave situation, it is the geodesic flow on $X$, while in
  the scattering case it is the Hamiltonian flow generated by the
  Hamiltonian $p(x,\xi)=\frac{|\xi|^2}{2}+V(x)$, in some energy interval
  $[E-\delta,E+\delta]$. Our assumptions on $X$ or $V(x)$
  imply that these classical flows are both ``strongly chaotic''.

In both cases, the long time properties of our quantum system
  involves a spectrum of complex ``eigenvalues'' associated with
  metastable states. 
Our main aim is to
understand the distribution of the lifetimes $\tau_j$ in the
semiclassical/high energy limit, especially the ones which are not
infinitesimally small when $\hbar\to 0$: we will thus focus on resonances
such that $|\Im z_j(\hbar)|=\mathcal{O}(\hbar)$.

\section{Transformation to nonselfadjoint spectral problems}
Each of these two quantum systems can be recast into a spectral
problem for an associated nonselfadjoint differential operator on $L^2$, with discrete spectrum
near the real axis.

\subsection{Damped quantum mechanics}\label{s:damped1}
Following \cite{Sjo00}, let us start from the damped wave system. 
The generalized eigenvalue equation
\eqref{e:k_j} for $\Re k_j\gg 1$ can be rewritten using an
effective ``Planck's constant'' $\hbar\approx (\Re k_j)^{-1}$, and
replacing $k_j$ by the ``energy'' $z_j= \frac{(\hbar k_j)^2}{2}=1/2+\mathcal{O}(\hbar)$:
\begin{equation}\label{e:P_b}
P_{dw}(\hbar) \psi_j = z_j\psi_j +\mathcal{O}(\hbar^2), \qquad P_{dw}(\hbar) = -\frac{\hbar^2\Delta}{2}- i \hbar b(x)\,.
\end{equation}
The principal symbol of the operator $P_{dw}(\hbar)$,
$p_0(x,\xi)=|\xi|^2/2$, is real and generates the
geodesic flow on $X$. The skew-adjointness of $P_{dw}(\hbar)$ only appears in the
subprincipal symbol $-i\hbar b(x)$: the latter does not influence the classical
dynamics, but is responsible for the decay of probability along the
flow. Indeed, for $\psi_0$ a wavepacket microlocalized
on $\rho_0=(x_0,\xi_0)\in T^*X$, its evolution
$\psi(t)=e^{-itP_{dw}(\hbar)/\hbar}\psi_0$ will be another
wavepacket microlocalized at $\rho_t=(x_t,\xi_t)=\Phi^t(\rho_0)$, with total probability
reduced by a finite factor 
\begin{equation}\label{e:damping-factor}
\frac{\|\psi(t)\|^2}{\|\psi_0\|^2}\approx \exp\Big(-2\int_0^t b(x_s)\,ds\Big)\,.
\end{equation}
In the limit $\hbar\to 0$, one can speak of a {\it damped classical dynamics}: each point $\rho_t$
evolves according to the geodesic flow, and carries a weight
which  gets reduced by the above factor along the flow. 

The horizontal spectral density of $P_{dw}(\hbar)$ is given by
Weyl's law, which (to lowest order) does not depend on the damping
\cite{Sjo00}: for any $c>0$,
\begin{equation}\label{e:Weyl}
\sharp\{{\rm Spec} (P_{dw}(\hbar))\cap ( [1/2-c \hbar,1/2+c\hbar] +i\mathbb{R})\} =
\hbar^{-d+1} \big(c\,C_X +\mathcal{O}(1)\big)\,,
\end{equation}
where $C_X>0$ only depends on $X$.
On the other hand, we will see that the distribution of the imaginary
parts $\Im z_j(\hbar)$
strongly depends on the {\it interplay} between the geodesic
flow and the damping $b(x)$.

\subsection{Complex scaled scattering Hamiltonian and open dynamics}
The scattering operator $P(\hbar)$ in \eqref{e:Schro-t} can be
transformed to a nonselfadjoint one through the complex scaling
method \cite{SjZw91}. One deforms the configuration space
$\mathbb{R}^d$ into a complex contour
$\Gamma_\theta=\{ x + i\theta f(x)\}$, where $f(x)=0$ for $x$ in a
ball $B(0,R)$ containing the support of $V(x)$ (the ``interaction region''), while
$f(x)=x$ for $|x|\geq 2R$. We take an angle $\theta= M
\hbar\log(1/\hbar)$, $M>0$ fixed. This
leads to a ``scaled'' operator $P_\theta(\hbar)$, which is
no more selfadjoint: in the sector
$-2\theta< \arg(z)\leq 0$ it admits discrete eigenvalues, which
correspond to the resonances of $P(\hbar)$. 

Our quest for resonances has turned into the spectral
study of $P_\theta(\hbar)$. This operator admits the symbol
$$
p_\theta(x,\xi) = p(x,\xi) - i \theta\langle\xi, df (x)^t \xi\rangle +
\mathcal{O}(\theta^2)|\xi|^2\,,\qquad p(x,\xi)=\frac{|\xi|^2}{2}+V(x) \,.
$$
The operator $P_\theta(\hbar)$
presents similarities with \eqref{e:P_b}: its principal symbol
$p(x,\xi)$ is real and
generates the Hamiltonian flow on $p^{-1}(E)$, while its
imaginary part is of higher order  $\hbar\log(1/\hbar)$, and generates a
damping outside $B(0,R)$. One
difference with $P_{dw}$ lies in the strength of this damping: for an initial wavepacket
localized on a point $\rho_0\in p^{-1}(E)$,
the full probability after time $t$ will be reduced by a factor
$\hbar^{2M\int_0^t\langle\xi_s,df(x_s)^t\xi_s\rangle ds}$: the probability is
semiclassically strongly suppressed as soon as the trajectory enters  the 
zone where $\langle\xi, df (x)^t \xi\rangle>0$ (this ``absorbing zone'' contains the
exterior of $B(0,2R)$). 

Classically, this corresponds to an ``open
dynamics'': the point $\rho_t$ evolves according to the Hamiltonian
flow, but it gets ``absorbed'', or ``killed'' as soon as it enters the 
absorbing zone. 

\subsubsection*{Trapped set}
For any energy $E>0$, the forward
(resp. backward)  trapped set $K_E^-$ (resp. $K_E^+$) is defined as
the set of initial points which remain bounded for all positive
(resp. negative) times:
$$
K_E^{\mp} = \{\rho_0=(x_0,\xi_0)\in p^{-1}(E),\ |x_t|\leq R,\ \forall t\gtrless 0\}\,,
$$
while the trapped set is made of their intersection
$K_E=K_E^-\cap K_E^+$. Notice that $K_E$ is a compact flow-invariant
set. The above remark shows that any wavepacket
localized on a point $\rho_0\not\in K_E^-$ will be absorbed after a finite
time, namely the time it takes to enter the absorbing zone. 
On the other hand, a point $\rho\in K_E^-$ will converge to the
trapped set $K_E$ as $t\to\infty$. This argument shows that, in some sense,
long time quantum mechanics (at energy $\approx E$) takes place
on $K_E$.

\section{Fractal Weyl laws}
The above argument can be made precise when estimating
the number of resonances of $P(\hbar)$ near $E$. In the case we are
interested in, $K_E$ is
a hyperbolic repeller (that is, there is no fixed point on $K_E$, and all trajectories are
hyperbolic), this number is bounded
above by a {\it fractal Weyl law} directly related with the geometry of
$K_E$ \cite{SjoZwo07}.
\begin{theorem} Assume that the trapped set $K_E$
  at some energy $E>0$ is a hyperbolic repeller, and write its Minkowski dimension
  ${\rm dim}_M(K_E) = 1+2\nu$.  
Then for any $c,\alpha>0$, one has in the semiclassical limit
\begin{equation}\label{e:fractalW}
\sharp \{{\rm Res}(P(\hbar))\cap ([E - c\hbar,E + c\hbar] - i\hbar[0,\alpha])\} =
\mathcal{O}(\hbar^{-\nu-0})\,.
\end{equation}
\end{theorem}
This theorem is proved by conjugating $P_\theta(\hbar)$ by a suitable
``weight'' $G(\hbar)$, so that the symbol $p_{\theta,G}$ of the
conjugated operator $e^{-G(\hbar)}P_\theta e^{G(\hbar)}$ satisfies
$\Im p_{\theta,G} < - 2C\hbar$ outside the
$\sqrt{\hbar}$-neighbourhood of $K_E$. Estimating the volume
of this neighbourhood, and some involved pseudodifferential calculus
on $p_{\theta,G}$, lead to
the above upper bound.

\medskip

A similar argument can be used to study the distribution of decay
rates $\Im z_j/\hbar$ for the operator \eqref{e:P_b}. For any
time $T>0$, one can construct a weight $G_T(\hbar)$, such that the conjugate
operator $P_{dw,G_T}(\hbar)=e^{-G_T(\hbar)}P_{dw}(\hbar)e^{G_T(\hbar)}$ admits the symbol
\begin{equation}\label{e:p_b,GT}
p_{dw,G_T}(x,\xi)=\frac{|\xi|^2}{2}-i\hbar b_T(x,\xi) +\mathcal{O}(\hbar^2)\,,
\end{equation}
where the subprincipal
symbol $b_T(\rho)= T^{-1}\int_{-T/2}^{T/2}b(x_t)\,dt$ is the
average of the damping along the flow. The geometric assumption of
negative curvature implies that the geodesic flow on $X$ is {\it
  Anosov}, in particular it is ergodic. This implies that, on the
energy shell $p^{-1}(1/2)$, the time
average $b_T(\rho)$ converges almost everywhere to the
microcanonical average $\bar b={\rm Vol}(X)^{-1}\int_X b(x)dx $ when
$T\to\infty$. From there one
can deduce that most of the eigenvalues
$\Re z_i(\hbar)\in [1/2,c\hbar,1/2+c\hbar]$ concentrate near
the ``typical line'' $\Im z = -\hbar\bar b$ \cite{Sjo00}.

Using finer properties of the Anosov geodesic flow, 
one can estimate the number of eigenvalues away
from this ``typical line'' \cite{Anan09}. To state the result, we need
to introduce the extremal ergodic averages of the damping, $b_{-}=\lim_{T\to\infty}\min_{p^{-1}(1/2)}b_T$,
and similarly for $b_+$. 
\begin{theorem}
Assume $X$ is a compact surface of negative curvature. Then, there exists a function
$H:\mathbb{R}\to \mathbb{R}$, strictly concave on $[b_-,
b_+]$ and equal to $-\infty$ outside,
with maximum $H(\bar b)=d-1$, such that for any
$c>0$ and any $\alpha\in [0,\bar b]$,
\begin{equation}\label{e:Fractal}
\sharp\{{\rm Spec} (P_{dw}(\hbar))\cap 
([1/2 - c\hbar,1/2 + c\hbar] - i\hbar[0,\alpha])\} = \mathcal{O}(\hbar^{-H(\alpha)-0})\,.
\end{equation}
A similar estimate holds for the range $-\Im z\geq
\hbar\alpha$, $\alpha > \bar{b}$.
\end{theorem}
Comparing this ``fractal Weyl
upper bound'' with
the Weyl law \eqref{e:Weyl}
confirms that most resonances are on the typical line.
The above theorem is obtained by studying the {\it large deviations} of the
value distribution of $b_T$ in the limit
$T\to\infty$: roughly speaking, for $\alpha <\bar b$, the volume of the points $\rho\in p^{-1}(1/2)$ such that
$b_T(\rho)\leq \alpha$ decays like $e^{T(H(\alpha)-(d-1))}$. This
volume estimate is then used to get \eqref{e:Fractal}.

In both situations, the upper bounds on counting resonances/eigenvalues
were obtained by deforming
the operator $P_\theta$ (resp. $P_{dw}$) by an appropriate microlocal weight, and
studying the imaginary part of the resulting operator by phase space
volume arguments, where the classical dynamics plays a prominent r\^ole.
 
\subsection{Spectral gaps}
We now present a complementary type of spectral information, which can be
obtained by a similar method for these two systems. Namely, we want to
understand if the lifetimes $\tau_j(\hbar)$ can be arbitrarily large in
the semiclassical limit; or on the opposite, if there exists a {\it gap} of
size $\propto\hbar$ between the real axis and the
eigenvalues/resonances. The presence of such a gap has important consequences on the
long time properties of the system.

Let us start with the scattering problem. The gap question can be
rephrased as: ``Are the metastable states able to concentrate on
$K_E$ when $\hbar\to 0$?'' The answer will result from
a {\it competition} between, one one side, the fast dispersion of wavepackets
due to the hyperbolic classical flow, on the other side
the ``thickness'' of the trapped set allowing the state to reconstruct
itself through constructive interferences.

A dynamical quantity reflecting this competition is of statistical
nature, it is a {\it topological pressure} associated with the flow on $K_E$:
$$
\mathcal{P}(-\varphi_u/2)=\lim_{T\to\infty} \frac1T
\log\sum_{\gamma:T\leq T_\gamma\leq T+1} \exp\Big(-\int_0^{T_\gamma} \varphi_u(\rho_t)/2\,dt\Big)\,.
$$
Each $\gamma$ is a periodic orbit in $K_E$ with period $T_\gamma$, and
$\varphi_u(\rho)$ is the unstable Jacobian of the flow.
The above mentioned competition lies in the fact that each
exponential weight gets very small when $T_\gamma\to\infty$, while the
number of terms grows exponentially.

The following gap criterion was first obtained \cite{Ikawa88} for the case of hard
obstacles, and then generalized in \cite{NoZw09} to smooth potentials.
\begin{theorem}
Assume the trapped set $K_E$ is a hyperbolic repeller. If the pressure
$\mathcal{P}(-\varphi_u/2)<0$, then for any $c,\epsilon>0$ and any
small enough $\hbar$, the strip \\$[E-c\hbar,E+c\hbar]+i\hbar
[\mathcal{P}(-\varphi_u/2)+\epsilon,0]$ does not contain any
resonance of $P(\hbar)$.
\end{theorem}
In dimension $d=2$, the condition
$\mathcal{P}(-\varphi_u/2)<0$ is equivalent with a purely
geometric statement, namely the fact that the Hausdorff dimension
${\rm dim}_H(K_E)<2$ (notice that ${\dim}_Hp^{-1}(E)=3$).

In the case of damped waves, the gap question
is nontrivial if $b_-=0$, that is, if there exists a flow-invariant subset
of $p^{-1}(1/2)$ with no damping. A result similar to the one above
was obtained in \cite{Schenck10}. In this case, the local decay of
probability is due to both hyperbolic dispersion
and damping.
\begin{theorem}
Assume that $X$ has negative curvature, and that the topological pressure
$\mathcal{P}(-\varphi_u/2 - b)<0$. Then, for any $c,\epsilon>0$,
and $\hbar>0$ small enough, the
strip $[1/2-c\hbar,1/2+c\hbar]+i\hbar[\mathcal{P}(-b-\varphi_u/2)+\epsilon,0]$
does not contain eigenvalues of $P_{dw}(\hbar)$.
\end{theorem}
\section{Open questions}
Most of the above results are upper or lower bounds. The natural
question is: ``Are these bounds sharp?''
The fractal Weyl bound \eqref{e:fractalW} is conjectured to be
sharp for $\alpha > 0$ large enough, a fact which has been tested numerically on a number of
examples, but could be proved only for a very specific toy model
\cite{NZ1}. 
On the opposite, the bounds \eqref{e:Fractal} for eigenmodes of the damped
wave equation are not
expected to be sharp for all values of $\alpha$.
The size of the gap itself is believed to be larger than the
topological pressure bound we gave above. Such an ``extra gap'' was
proved for the 3-disk scattering, using advanced estimates on
classical mixing \cite{PetStoy10}.

\smallskip

\noindent
{\sc Acknowledgments.}
The author has been partially supported by the Agence Nationale
de la Recherche under the grant ANR-09-JCJC-0099-01.

\end{document}